# Decision Approach and Empirical Bayes FCR-Controlling Interval for Mixed Prior Model

## Zhigen Zhao

*e-mail:* `zhigen.zhao@gmail.com`

**Abstract:** In this paper, I apply the decision theory and empirical Bayesian approach to construct confidence intervals for selected populations when true parameters follow a mixture prior distribution. A loss function with two tuning parameters $k_1$ and $k_2$ is coined to address the mixture prior. One specific choice of $k_2$ can lead to the procedure in Qiu and Hwang (2007); the other choice of $k_2$ provides an interval construction which controls the Bayes FCR. Both the analytical and extensive numerical simulation studies demonstrate that the new empirical Bayesian FCR controlling approach enjoys great length reduction. At the end, I apply different methods to a microarray data set. It turns out that the average length of the new approach is only 57% of that of Qiu and Hwang's procedure which controls the simultaneous non-coverage probability and 66% of that of Benjamini and Yekutieli (2005)'s procedure which controls the frequentist's FCR.

**AMS 2000 subject classifications:** Decision Bayes, Loss Function, Simultaneous Intervals..

## 1. Introduction

Simultaneous interval estimation for a large number of selected parameters is challenging especially when the number of observations for each parameter is very small. The difficulties are the selection bias (see Qiu and Hwang 2007 and Hwang 1993) and the multiplicity. The traditional approach, which treats all the parameters as fixed, seems to have little power when the dimension tends to be very large, for instance, several thousands in microarray. However, the empirical Bayesian approach is known to be able to *borrow strength* across the populations. Thus, it is very likely that this method will provide us with some satisfactory procedures.

In the past, people attempted to estimate the parameters for selection populations (see, for example, Cohen and Sackrowitz 1982 and Hwang 1993). However, very few people knew how to construct an interval for selected population. The first exciting work was written by Benjamini and Yekutieli (2005) (I will use B-Y (2005) to represent this work throughout the paper). They adapted the concept of FDR from multiple testing and coined a concept False Coverage Rate (FCR) for simultaneous intervals. This criterion is much less conservative than the simultaneous non-coverage coefficient. They have constructed confidence intervals for multiple selected parameters which can control the FCR at specified $q$-level, typically 5%. They centered their intervals upon the estimator $X_i$'s which are biased for selected populations and addressed the multiplicity by lengthening the intervals. Consequently, their intervals have extremely large average half length.

---

*This is an original survey paper

In 2008, Zhao and Hwang introduced the Bayes FCR and connected Bayes confidence interval which aims at controlling Bayesian non-coverage coefficients with the Bayesian FCR controlling intervals. They applied this general theorem to the normal-normal setting where the observations follow a normal distribution with unequal but know variances and the parameters follow a normal prior. They used the empirical Bayesian approach to derive explicit intervals which can control the empirical Bayes FCR. Their construction reduced the average length of B-Y's procedure dramatically because they addressed the multiplicity by modifying the centers instead of the lengths.

Another exciting work is Qiu and Hwang (2007), which offers a way to construct intervals that can control the simultaneous coverage coefficient for selected popultions. Other than the normal-normal model, they treated the so-called normal-mixture model where the prior distribution of the true parameters is a mixture of a normal random variable with an equal, known variance and a single point *zero*. Because they have addressed the multiplicity by Bonferroni's correction, their lengths tend to be large when many parameters are selected.

In this paper, I use the decision approach and empirical Bayes to construct intervals for selected populations under the same model setting of Qiu and Hwang (2007). Application of decision approach to interval/set estimation has a long history which dates back to Faith (1976), Casella and Hwang (1983), and He (1992). Recently, Hwang, Qiu, and Zhao (2008) have constructed the double shrinkage empirical confidence interval for one single parameter when assuming the variances to be unequal and unknown. However, all the loss functions they have used are not appropriate under the mixed prior model (Detailed argument is in section 2.2). Thus a new loss function with two tuning parameters $k_1$ and $k_2$ is proposed and strongly recommended. One specific choice of $k_2$ results in Qiu and Hwang (2007)'s procedure. The other choice of $k_2$ provides us with a way to construct the empirical Bayesian FCR-controlling intervals based on the normal-mixture model.

In section 2, I introduce the model setting and the decision Bayes rule based on our new loss function. In section 3, I will connect the decision Bayesian rule with Qiu and Hwang (2007)'s procedure first and then derive a procedure which can control the Bayes FCR. In section 4, empirical Bayesian approach is constructed and evaluated both numerically and analytically. In section 5, I apply the confidence intervals constructed in section 4 to a real microarray dataset and compare it with B-Y (2005)'s and Qiu and Hwang's procedures. It turns out that my procedure out-performs theirs. The average length of my interval is only 57% of that of Qiu and Hwang's (2007) procedure which controls the simultaneous coverage probability and 66% of that of B-Y (2005)'s procedure which controls the frequentist's FCR.

## 2. Normal-Mixture Model for the means

### 2.1. Model Assumption

In microarray, it is generally assumed that observed differentially expressed levels $X_i$'s are normally distributed with true means $\theta_i$'s, $i = 1, 2, \cdots, p$, where the dimension $p$ varies from several thousands to over thirty thousands. Due to the extremely large number of dimensions, it is strongly recommended to use a prior to model the true means $\theta_i$'s. A natural choice is the normal prior where $\theta_i \overset{i.i.d.}{\sim} N(0, \tau^2)$.

However, in Qiu and Hwang (2007), they have applied the Q-Q plot to a microarray data and shown that normal-normal model cannot fit the data well. To remedy this, they introduced the *normal-mixture model* as following. Assume that $X_i|\theta_i \sim N(\theta_i, \sigma^2)$, and

$$\pi(\theta_i) \begin{cases} = 0 & \text{with probability } \pi_0, \\ \sim N(0, \tau^2) & \text{with probability } \pi_1 = 1 - \pi_0. \end{cases} \quad (1)$$

I use an indicator function $I_i$ to describe whether $\theta_i$ is 0, i.e. $I_i = 0$ if $\theta_i = 0$ and $I_i = 1$ if $\theta_i \sim N(0, \tau^2)$. Initially, I assume that hyper parameters $\tau^2$ and $\pi_0$ are known and derive the corresponding decision Bayesian procedure. In section 4, I estimate them through data by using consistent estimators and derive a empirical Bayesian procedure.

### 2.2. Bayes Interval

In history, there are many attempts to apply the decision Bayes approach to construct confidence sets/intervals. Faith (1976) first introduced a linear loss function for confidence set $CI$ of the parameter $\theta$ as $L(\theta, CI) = kVolume(CI) - I_{CI}(\theta)$ where the tuning parameter $k$ was determined by some minimax rule in Casella and Hwang (1983). He (1992) used $L_i(\theta_i, CI_i) = kLen(CI_i) - I_{CI_i}(\theta_i)$ as the loss function for the interval estimator $CI_i$ of the parameter $\theta_i$. Hwang Qiu and Zhao (2008) modified the loss function above as $L(\theta_i, CI_i) = \frac{k}{\sigma_i}Len(CI_i) - I_{CI_i}(\theta_i)$ and constructed the double shrinkage confidence interval when assuming variances to be unequal and unknown. However, all these loss functions are not appropriate for normal-mixture model (1). In fact, for any given confidence interval, one can construct a new interval, which is the union of the existing procedure and *zero*. This new approach boosts the coverage probability while causes no change of the length. Consequently, the conditional expected loss of the new construction is always less or equal than that of the original approach. As a result, the decision Bayes suggests that *zero* should be included in every interval. But practically, such constructions have no power and appear to be useless.

In order to avoid this phenomenon, I put extra terms which influence the loss function only when the point *zero* is included and thus define the loss function as,

$$L(\theta_i, CI_i) = k_1 Len(CI_i)I(I_i = 1) - I_{CI_i}(\theta_i, I_i = 1) + I_{CI_i}(0)(k_2 - I(I_i = 0)), 0 \leq k_2 \leq 1. \tag{2}$$

The first two terms balance the length and the true coverage when the true parameters $\theta_i$'s are generated from the normal random variable. The tuning parameter $k_1$ will be determined later in this section. The last two terms affect the loss function only when the corresponding interval does include *zero*. Upon this, if $\theta_i$ is indeed *zero*, then $k_2 - I(I_i = 0) = k_2 - 1 \leq 0$, implying that including *zero* is useful. On the other hand, if $\theta_i$ is not *zero*, then $k_2 - I(I_i = 0)$ is positive and becomes a penalty term. Thus, appropriate choice of the tuning parameter $k_2$ guides us to decide when *zero* should be included.

Furthermore, the flexibility of choosing $k_2$ offers us constructions under different settings. For example, when assuming the normal-normal model, the loss function (2) reduces to He (1992)'s if I set $k_2 = 0$. In section 3, I apply two different choices of $k_2$, one of which will reproduce Qiu and Hwang (2007)'s procedure and the other of which provides a construction that can control the Bayesian FCR at $q$-level.

Now, I have all the pieces to construct the decision Bayes rule, i.e. I want to construct a Bayes interval $CI_i^B$ such that it minimizes $E(L(\theta_i, CI_i|X))$ for any observation $X$ when assuming the normal-mixture model (1) and the loss function (2).

**Theorem 2.1** *Let $\pi_i^0(X) = P(\theta_i = 0|X) = P(I_i = 0|X)$ and $\pi_i^1(X) = 1 - \pi_i^0(X)$. Then*

$$EL(\theta_i, CI_i|X) = \pi_i^1(X) \int_{CI_i} (k_1 - \pi(\theta_i|X, I_i = 1))d\theta_i + I_{CI_i}(0|X)(k_2 - \pi_i^0(X)). \tag{3}$$

*The Bayes interval is*

$$CI_i = \begin{cases} \{\theta_i : k_1 < \pi(\theta_i|X_i, I_i = 1)\} \setminus \{0\} & \text{if } k_2 > \pi_i^0(X), \\ \{\theta_i : k_1 < \pi(\theta_i|X_i, I_i = 1)\} \cup \{0\} & \text{if } k_2 \leq \pi_i^0(X). \end{cases} \tag{4}$$

Intuitively, for any given observation $X_i$, if the conditional probability $P(\theta_i = 0|X)$ is small, it is very unlikely that $\theta_i = 0$ and *zero* should be included. On the other hand, larger $\pi_i^0(X)$ indicates that *zero* should be included. Theorem 2.1 tells us that the parameter $k_2$ is the threshold value.

Under model (1), $\pi(\theta_i|X, I_i = 1) \sim N(MX_i, M\sigma^2)$ where $M = \frac{\tau^2}{\tau^2 + \sigma^2}$, therefore

$$\{\theta_i : k_1 < \pi(\theta_i|X_i, I_i = 1)\} = \{\theta_i : (\theta_i - MX_i)^2 < -M\sigma^2(2\log k_1\sqrt{2\pi} + \log M\sigma^2)\}.$$

As in the Section 3 of Hwang, Qiu, and Zhao (2008), one wants to obtain a traditional normal interval when the non-informative prior is applied, i.e., if setting $\tau \to \infty$, $M \to 1$, one wants the corresponding interval $\{\theta_i : \frac{(\theta_i - X_i)^2}{\sigma^2} < -(2\log k_1\sqrt{2\pi} + \log \sigma^2)\}$ to coincide with normal interval $(X_i - z_{q/2}\sigma, X_i + z_{q/2}\sigma)$ where $z_{q/2}$ is the critical value such that $P(|Z| > z_{q/2}) = q$ when $Z$ is a standard normal random variable. Therefore, the constant $k_1$ should be chosen such that $z_{q/2}^2 = -(2 \log k_1 \sqrt{2\pi} + \log \sigma^2)$. Plug this constant $k_1$ back to Bayes interval (4). Then the decision Bayes interval becomes

$$CI_i^B = \begin{cases} \{\theta_i : (\theta_i - MX_i)^2 < M\sigma^2(z_{q/2}^2 - \log M)\} \setminus \{0\} & \text{if } k_2 > \pi_i^0(X), \\ \{\theta_i : (\theta_i - MX_i)^2 < M\sigma^2(z_{q/2}^2 - \log M)\} \cup \{0\} & \text{if } k_2 \le \pi_i^0(X). \end{cases} \quad (5)$$

Unlike the interval $MX_i \pm \sqrt{M}\sigma z_{q/2}$, which is directly derived from the posterior distribution, the major part of (5) has an extra positive term $M\sigma^2(-\log M)$ which is necessary to boost the coverage probability when the hyper parameters are estimated through the data in section 4. In the next section, I will choose the value of the parameter $k_2$ under two different problem settings and derive the decision Bayes interval accordingly.

## 3. Choose $k_2$

### 3.1. *Qiu and Hwang (2007)*

Qiu and Hwang (2007) constructed the interval for $K$ selected populations $X_{(p-K+1)}, X_{(p-K+2)}, \cdots, X_{(p)}$ under the model (1) where the observations $X_{(j)}$'s satisfy

$$|X_{(1)}| \le |X_{(2)}| \le \cdots \le |X_{(p)}|.$$

Assume $\theta_{(j)}$ is the true parameter that corresponding to the observation $X_{(j)}$. Note that $|\theta_{(p)}|$ is not necessarily equal to $\max\limits_{1 \le j \le p} |\theta_{(j)}|$. I construct the interval for $\theta_{(j)}$ where $p - K + 1 \le j \le p$ as

$$CI_{(j)}^B = \begin{cases} \{\theta_{(j)} : (\theta_{(j)} - MX_{(j)})^2 < M\sigma^2(z_{q/2K}^2 - \log M)\} \setminus \{0\} & \text{if } k_2 > \pi_{(j)}^0(X), \\ \{\theta_{(j)} : (\theta_{(j)} - MX_{(j)})^2 < M\sigma^2(z_{q/2K}^2 - \log M)\} \cup \{0\} & \text{if } k_2 \le \pi_{(j)}^0(X). \end{cases} \quad (6)$$

When compared with (5), the major difference is that I use the critical value $z_{q/2K}$ to address the multiplicity. This is known as Bonferroni's correction.

Direct calculation shows that for each $j$,

$$P(\theta_{(j)} \notin CI_{(j)}^B | X) \le q/K + \pi_{(j)}^0(X)(I(\pi_{(j)}^0(X) < k_2) - q/K).$$

Consequently, the simultaneous non-coverage coefficient satisfies

$$P(\theta_{(j)} \notin CI_{(j)}, j = p-K+1, \cdots, p | X) \le q + \sum_{j=p-K+1}^{p} \pi_{(j)}^0(X)(I(\pi_{(j)}^0(X) < k_2) - q/K). \quad (7)$$

If $k_2$ is chosen to be the maximum $k$ such that the summation above is non-positive, i.e.

$$k_2 = \max_k \{ \sum_{j=p-K+1}^{p} \pi_{(j)}^0(X)(I(\pi_{(j)}^0(X) < k) - q/K) \le 0 \}. \quad (8)$$

Then the non-coverage coefficient $P(\theta_{(j)} \notin CI_{(j)}, j = p - K + 1, \cdots, p)$ is controlled at the $q$-level. Surprisingly, this choice of $k_2$ is exactly the same as in Theorem 4 of Qiu and Hwang (2007). Therefore, Qiu and Hwang (2007)'s Bayes procedure is exactly the same as (6).

### *3.2. Bayes FCR Controlling Interval*

Benjamini and Yekutieli (2005) initiated the concept of FCR, which is much less conservative than the simultaneous non-coverage coefficients. Zhao and Hwang (2008) have extended this idea to the Bayesian framework through a new concept, Bayes FCR. They have shown that there is a natural connection between the Bayes FCR and the Bayes confidence interval. In this subsection, I will show that (5) can control the Bayes FCR at the $q$-level if $k_2$ is chosen appropriately.

**Theorem 3.1** *Assume that $\mathcal{R}(X)$ is the index set of observations that are selected for interval estimation. $R = \#\mathcal{R}$. Define*

$$f(p, \tau^2, \pi_0, k) = E(\sum_{i \in \mathcal{R}} \frac{\pi_i^0(X)(I(\pi_i^0(X) < k) - q)}{R} I(R > 0)),$$

*and $k_2 = \max\limits_k \{k, f(p, \tau^2, \pi_0, k) \leq 0\}$. Then intervals (5) satisfies*

$$FCR_\pi \leq qP(R > 0).$$

*In other words, the Bayes FCR of the intervals (5) has been controlled at $q$ level.*

Now assume that the selection rule in Qiu and Hwang (2007) is applied, i.e., the last $K$ observations after ordering $X_1, X_2, \cdots, X_p$ according to their absolute value increasingly are selected. Then $f(p, \tau^2, \pi_0, k) = E(\sum_{i=p-K+1}^{p} \frac{\pi_{(i)}^0(X)(I(\pi_{(i)}^0(X) < k) - q)}{K})$. Comparing it with the expectation part in (8), $f$ is always smaller when $K > 1$, which implies that the choice of $k_2$ for controlling Bayes FCR interval is always larger than the choice of Qiu and Hwang's. Consequently, under the same setting, the frequency that (5) includes *zero* is less than that of Qiu and Hwang's. Furthermore, since they addressed the multiplicity by Bonferroni's correction, the half length $M\sigma^2(z_{q/2K} - \log M)$ is much larger than the half length of Bayes FCR controlling interval (5) and the discrepancy becomes large when $K$ is big. These two facts all implied that Bayes FCR controlling interval is less conservative than Qiu and Hwang (2007).

Another nice thing about this theorem is that it holds for any selection rule, including pre-determined and data-driven selection rule. For example, when observations are selected according to Benjamini and Hochberg (1995)'s procedure which controls the False Discovery Rate at $q$-level and $k_2$ is simulated accordingly, the above theorem still guarantees that (5) controls the Bayes FCR at $q$-level.

The choice of $k_2$ depends on the unknown expectation, which prevents us from finding $k_2$ explicitly. However, $k_2$ can be easily determined by simulation once the hyper-parameters are known.

## 4. Empirical Bayes Approach

In this section, I estimate unknown hyper-parameters through the data and obtain a practical confidence interval. Our goal is to construct the confidence intervals for selected parameters such that the Bayes FCR can always be controlled for a class of prior distributions which are determined by the hyper-parameters $\pi_0$ and $\tau^2$. This approach is named *empirical Bayes FCR controlling intervals*, according to Zhao and Hwang (2008).

Recall the model 1. Then $EX_i^2 = \sigma^2 + \pi_1\tau^2$, and $EX_i^4 = 3(\sigma^4 + 2\pi_1\sigma^2\tau^2 + \pi_1\tau^4)$. By using the method of moments, one could get reliable estimators of $\pi_0$ and $\tau^2$ when $p$ is sufficiently large,

$$\hat{\pi}_1 = \frac{(m_2 - \sigma^2)^2}{m_4/3 + \sigma^4 - 2\sigma^2 m_2}, \hat{\tau}^2 = \frac{(m_2 - \sigma^2)}{\hat{\pi}_1}. \qquad (9)$$

Plug these two estimators back to the function of $f$ and simulate the value of $k_2$, which is denoted by $\hat{k}_2$. Assume that $\hat{M}$ and $\hat{\pi}_i^0(X)$ are the estimators of $M$ and $\pi_i^0(X)$ when $\pi_0$ and $\tau^2$ are replaced by (9). Then I can construct the empirical Bayes interval as,

$$CI_i^{EB} = \begin{cases} \{\theta_i : (\theta_i - \hat{M}X_i)^2 < \hat{M}\sigma^2(z_{q/2}^2 - \log\hat{M})\} \setminus \{0\} & \text{if } \hat{k}_2 > \hat{\pi}_i^0(X), \\ \{\theta_i : (\theta_i - \hat{M}X_i)^2 < \hat{M}\sigma^2(z_{q/2}^2 - \log\hat{M})\} \cup \{0\} & \text{if } \hat{k}_2 \leq \hat{\pi}_i^0(X). \end{cases} \qquad (10)$$

The following theorem describes the asymptotic property of the construction.

**Theorem 4.1** *For any $0 < \pi_0 < 1$, $\tau^2 > 0$, if $\forall \epsilon > 0$, there $\exists \delta, N > 0$, such that $\forall p > N, k, k^{'} > 0$, $(\tau^{'2} - \tau^2)^2 + (\pi_0' - \pi_0)^2 + (k' - k)^2 < \delta$ implies*

$$|f(p, \tau^{'2}, \pi_0', k') - f(p, \tau^2, \pi_0, k)| < \epsilon. \qquad (11)$$

*Then under the model (1), the empirical Bayes interval (10) satisfies*

$$\limsup_{p \to \infty} FCR_\pi \leq q.$$

**Proposition 4.1** *If the selected parameters are the first $R$ parameters and $R \to \infty$ when $p \to \infty$. Then $f$ satisfies the condition (11).*

This proposition implies that when all observations are selected for interval estimation, (10) can control the empirical Bayes FCR asymptotically.

However, like all other existing constructions such as Casella and Hwang (1983), Qiu and Hwang (2007), and Hwang, Qiu, and Zhao(2008), the interval (10) cannot provide a satisfactory answer automatically for the finite sample case.

In figure 1, I have plotted a figure of Bayes FCR of the empirical Bayes interval versus the procedure of B-Y under different settings of hyper-parameter $(\pi_0, \tau^2)$ when $p = 1000$ and only the top 100 observations are selected for interval estimation. B-Y's procedure can always control the FCR at the 5% level; however, their procedures are way too conservative in terms of extremely low Bayes FCR when $M$ is close to 1 and large average length. The green line, corresponding to the construction (10), performs well when $\tau^2$ is relatively large; however some modifications are strongly required when $\tau^2$ is small.

Qiu and Hwang (2007) have argued that $\pi_0$ is nearly unidentifiable when $\tau$ is small. This will cause the estimator (9) to be very inaccurate. Therefore, they mixed their empirical Bayes intervals with the Bonferroni correction $(X_i - z_{q/2p}\sigma, X_i + z_{q/2p}\sigma)$ based on a threshold: $\min(\sqrt{720/p}, 0.6)$, obtained from extensive numerical calculations. It also seems necessary to mix the procedure (10) with the interval $(X_i - z_{Rq/2p}\sigma, X_i + z_{Rq/2p}\sigma)$, which is inspired by B-Y (2005). The following analytic argument can help us to find the threshold value much easier than Qiu and Hwang (2007).

Recall that $EX_i^2 = \sigma^2 + \pi_1\tau^2$ and $EX_i^4 = 3(\sigma^4 + 2\pi_1\sigma^2\tau^2 + \pi_1\tau^4)$, therefore $\tau^2 + 2\sigma^2 = \frac{EX_i^4/3 - \sigma^4}{EX_i^2 - \sigma^2}$. Use $m_2 = \sum X_i^2/p$ and $m_4 = \sum X_i^4/p$ to denote the second and fourth moments, then $\hat{\tau}^2 + 2\sigma^2 = \frac{m_4/3 - \sigma^4}{m_2 - \sigma^2}$.

Since the left hand side is always greater or equal than $2\sigma^2$, $\tau^2$ is not estimable when the right hand side is less than $2\sigma^2$. Therefore, I can carefully choose a proper $\tau_0^2$, such that the probability of the right hand side is smaller than $2\sigma^2$, i.e. the probability that $\pi_0$ and $\tau^2$ are not estimable, is controlled at the level of $q$. Therefore, set the threshold value $\tau_0^2$ to satisfy $P_{\tau^2 = \tau_0^2}\left(\frac{m_4/3 - \sigma^4}{m_2 - \sigma^2} \leq 2\sigma^2\right) \leq q$.

Now consider the special case when $\pi_1 = 1$ and calculate $\tau_0^2$. Use $m_4'$ and $m_2'$ to denote the second and fourth moments of the standard normal distribution when there are $p$ observations. Then $m_4 = (\tau^2 + \sigma^2)^2 m_4'$ and $m_2 = (\tau^2 + \sigma^2)m_2'$. I choose $\tau_0^2$ such that

$$P_{\tau^2 = \tau_0^2}\left((\tau^2 + \sigma^2)^2\frac{m_4'}{3} - 2\sigma^2(\tau^2 + \sigma^2)m_2' + \sigma^4 < 0\right) \leq q$$

by simulation.

Based on the cutoff, the final empirical Bayes FCR controlling interval with mixture is defined as

$$CI_i^{Final} = \begin{cases} X_i \pm z_{Rq/(2p)}\sigma & \text{if } m_2 - \sigma^2 < \tau_0^2, \\ CI_i^{EB}, & \text{if } m_2 - \sigma^2 > \tau_0^2. \end{cases} \tag{12}$$

In figure 1, the red solid line corresponds to the above empirical Bayes intervals. They perform the same as BY when $\tau^2$ is very small because of the mixed procedure. The portion of the mixture increases when $\pi_0$ increases. However, (12) performs better than theirs when $\tau^2$ is larger. The discrepancy is significant when $M \to 1$.

I have also plotted the simulated average length in figure 2 that corresponds to the same model settings in figure 1. The average length of (12) is uniformly less or equal than the average length of B-Y's procedure. The ratio can even be as small as 56%.

In figures 3 and 4, I repeat the simulation setting but change the selection rule to Benjamini and Hochberg (1995)'s procedure which aims at finding significant observations while controlling the False Discovery Rate at 5%-level. The intervals (12) can control the empirical Bayesian FCR at 5%-level based on this data-driven selection. Compared with B-Y (2005)'s procedure, the improvement of the average length is even more significant than that corresponding to the fixed selection rule. The ratio can be as small as 43%.

## 5. Real Data Analysis

In this section, I apply different intervals to a microarray data set, the Synteni data of Kerr, Martin, and Churchill (2000), which was revisited by Hsu et al.(2006) and Qiu and Hwang (2007). The description of the data set can be found in Kerr, et al. (2000). The figure 6 of Qiu and Hwang (2007) is a Q-Q plot of the ANOVA estimator $X_g$, which shows that the normal-mixture model (1) fits the data well.

Hsu et al. (2006) uses simultaneous confidence intervals to detect genes with an expression level of $\Delta = 3$ or more. I will first apply Benjamini and Hochberg (1995)'s procedure to select parameters with expression levels significantly larger or equal than $\log_2 3$, and then construct the simultaneous interval for such selected observations. B-H's procedure declares that the first 89 genes are significant.

In figure 5, I construct the confidence intervals for these 89 genes by using Qiu and Hwang (2007)'s, B-Y's, (12). Our confidence interval (12) for $\theta_{(g)}$ is $0.93X_{(g)} \pm 0.96$. Compared with the interval $X_{(g)} \pm 1.47$ of BY's procedure, $0.93X_{(g)} \pm 1.67$ of Qiu and Hwang (2006), our intervals enjoy great length reduction.

## 6. Discussion

In this article, I have defined a new loss function for confidence interval construction when assuming the mixed prior model (1). I use two different ways to choose the tuning parameter in the loss function to obtain Qiu and Hwang (2007)'s procedure and the empirical Bayesian FCR controlling intervals. I conclude that the new empirical Bayesian FCR controlling interval is better than other existing procedures because of the sharp improvement over the average length.

However, there are still much need for further research. In model (1), I assume equal and known variance $\sigma^2$, which is not generally a practical assumption. Hwang, Qiu, and Zhao (2008) proposed a double shrinkage empirical Bayesian interval for single parameter without selection under the normal-lognormal model. Therefore, one natural extension of this work is to consider the mixture-prior model when variances are unequal and unknown. The loss function (2) provides us a potential tool to construct corresponding intervals.

**Supplementary Materials**

*A.1. Technical Details of Mathematical Results*

**Proof of Theorem 2.1.**

Firstly,

$$EL(\theta_i, CI_i|X) \tag{A.1}$$

$$= k_1 Len(CI_i)P(I_i = 1|X) - \int I(\theta_i \in CI_i, I_i = 1)m(\theta_i|X)d\theta_i + I_{CI_i}(0|X)(k_2 - \pi_i^0(X)).$$

The integration $\int I_{CI_i}(\theta_i, I_i = 1)m(\theta_i|X)d\theta_i$ can be written as $\int_{CI_i} m(\theta_i, I_i = 1|X)d\theta_i$ where $\pi(\theta_i, I_i = 1|X) = \pi_i^1(X)\pi(\theta_i|I_i = 1, X)$. Write $Len(CI_i)$ as $\int_{CI_i} 1d\theta_i$. Then (A.1) equals to

$$\pi_i^1(X)\int_{CI_i}(k_1 - \pi(\theta_i|X, I_i = 1))d\theta_i + I_{CI_i}(0|X)(k_2 - \pi_i^0(X)). \tag{A.2}$$

Now consider two intervals $CI_i^1$ and $CI_i^2$ where $CI_i^1 = \{\theta_i : k_1 < \pi(\theta_i|X, I_i = 1)\} \setminus \{0\}$ and $CI_i^2 = \{\theta_i : k_1 < \pi(\theta_i|X, I_i = 1)\} \cup \{0\}$. Then both $CI_i^1$ and $CI_i^2$ minimize the first term of the formula (A.2). Since $0 \in CI_i^2$ and $0 \notin CI_i^1$, then

$$EL(CI_i^2|X) = EL(CI_i^1|X) + (k_2 - \pi_i^0(X)).$$

Consequently, the Bayes interval includes 0 if and only if $k_2 < \pi_i^0(X)$, i.e. it is the one that is defined in (4).

**Proof of Theorem 3.1.**

According to Zhao and Hwang (2008),

$$FCR_\pi = E\frac{\sum_{i \in \mathcal{R}} P(\theta_i \notin CI_i|X)}{R}I(R > 0).$$

Since

$$P(\theta_i \notin CI_i^B|X)$$
$$= P(\theta_i \notin CI_i^B|X, I_i = 0)P(I_i = 0|X) + P(\theta_i \notin CI_i^B|X, I_i = 1)P(I_i = 1|X)$$
$$= \pi_i^0(X)I(\pi_i^0(X) < k_2) + (1 - \pi_i^0(X))P(\theta_i \notin CI_i^B|X, I_i = 1), \tag{A.3}$$

and $P(\theta_i \notin CI_i^B|X, I_i = 1) \le q$,

$$FCR_\pi \le qE(I(R > 0)) + E\frac{\sum_{i \in \mathcal{R}} \pi_i^0(X)(I(\pi_i^0(X) < k_2) - q)}{R}I(R > 0) = qP(R > 0) + f(k_2).$$

The choice of $k_2$ ensures that $f(k_2) \le 0$. Consequently,

$$FCR_\pi \le qP(R > 0).$$

**Proof of theorem 4.1.**

Before the proof, I will state and prove the following lemma.

**Lemma A.1** *Assume that $\hat{\tau}^2$ and $\hat{\pi}_0$ are consistent estimators of $\tau^2$ and $\pi_0$, then for any $\delta > 0$, there $\exists P_0 > 0$ such that $\forall p > P_0$,*

$$|\hat{\pi}_i^0 - \pi_i^0| \le \delta, \text{ for all } i = 1, 2, \cdots, p.$$

Direct calculation shows that $\pi_i^0 = \frac{\pi_0}{\pi_0 + \pi_1 \frac{\sigma}{\sqrt{\sigma^2 + \tau^2}} \exp(\frac{MX_i^2}{2\sigma^2})}$ and $\hat{\pi}_i^0$ has the same form as $\pi_i^0$ except that $\pi_0$ and $\tau^2$ are replaced by their estimators $\hat{\pi}_0$ and $\hat{\tau}^2$. Now, I introduce an intermediate estimator $\tilde{\pi}_i^0$ where $\pi_0$ is assumed known. I shall prove that the lemma holds for $\tilde{\pi}_i^0$ first.

Since $\hat{\tau}^2$ is consistent, $\hat{M} = \frac{\hat{\tau}^2}{\hat{\tau}^2 + \sigma^2}$ is also a consistent estimator of $M$. Then, for $\epsilon = \frac{1}{k} < \min(\frac{1-M}{M}\delta, \frac{\pi_1 \sigma}{\pi_0\sqrt{\sigma^2 + \tau^2}}\delta)$, there exists $N$, such that $\forall p > N$, $|\hat{M} - M| < \epsilon M$.

Without loss of generality, assume that $M > \hat{M}$, i.e. $0 < M - \hat{M} < \epsilon M = \frac{M}{k}$. Since $M$ is a increasing function with respect to $\tau^2$ when $\sigma^2$ is fixed, therefore $\tau^2 > \hat{\tau}^2$. Direct calculation shows that

$$\tilde{\pi}_i^0 - \pi_i^0 = \frac{\pi_0 \pi_1 \sigma (\sqrt{\frac{\sigma^2 + \hat{\tau}^2}{\sigma^2 + \tau^2}} \exp(\frac{(M - \hat{M})X_i^2}{2\sigma^2}) - 1)}{(\pi_0 \sqrt{\sigma^2 + \hat{\tau}^2} \exp(-\frac{\hat{M}X_i^2}{2\sigma^2}) + \pi_1 \sigma)(\pi_0 + \pi_1 \frac{\sigma}{\sqrt{\sigma^2 + \tau^2}} \exp(\frac{MX_i^2}{2\sigma^2}))} \quad \text{(A.4)}$$

Since $0 < \hat{M} < M$,

$$0 < \frac{\sigma^2 + \hat{\tau}^2}{\sigma^2 + \tau^2} = \frac{1 - M}{1 - \hat{M}} < 1.$$

Consequently,

$$\sqrt{\frac{\sigma^2 + \hat{\tau}^2}{\sigma^2 + \tau^2}} > \frac{\sigma^2 + \hat{\tau}^2}{\sigma^2 + \tau^2} = \frac{1 - M}{1 - \hat{M}}.$$

Therefore, (A.4) implies that

$$\tilde{\pi}_i^0 - \pi_i^0 > \frac{\pi_0 \pi_1 \sigma (\frac{1-M}{1-\hat{M}} - 1)}{(\pi_0 \sqrt{\sigma^2 + \hat{\tau}^2} \exp(-\frac{\hat{M}X_i^2}{2\sigma^2}) + \pi_1 \sigma)(\pi_0 + \pi_1 \frac{\sigma}{\sigma^2 + \tau^2} \exp(\frac{MX_i^2}{2\sigma^2}))}$$

Since the numerator is negative and the denominator is larger than $\pi_0 \pi_1 \sigma$,

$$\tilde{\pi}_i^0 - \pi_i^0 > \frac{\pi_0 \pi_1 \sigma \frac{\hat{M} - M}{1 - \hat{M}}}{\pi_0 \pi_1 \sigma} > \frac{\hat{M} - M}{1 - \hat{M}}.$$

Furthermore, $\hat{M} - M > -\epsilon M$ implies that

$$\tilde{\pi}_i^0 - \pi_i^0 > \frac{M}{1 - M}(-\epsilon) > -\delta. \quad \text{(A.5)}$$

On the other hand,

$$\tilde{\pi}_i^0 - \pi_i^0 \le \frac{\pi_0 \pi_1 \sigma (\exp(\frac{\epsilon M X_i^2}{2\sigma^2}) - 1)}{\frac{\pi_1^2 \sigma^2}{\sqrt{\sigma^2 + \tau^2}} \exp(\frac{MX_i^2}{2\sigma^2})} = \frac{\pi_0 \sqrt{\sigma^2 + \tau^2}}{\pi_1 \sigma} \cdot \frac{exp(\frac{\epsilon M X_i^2}{2\sigma^2}) - 1}{\exp(\frac{MX_i^2}{2\sigma^2})}.$$

I use $C$ to denote the constant $\frac{\pi_0\sqrt{\sigma^2+\tau^2}}{\pi_1\sigma}$, and let $y = \exp(\frac{\epsilon MX_i^2}{2\sigma^2})$, then $\exp(\frac{MX_i^2}{2\sigma^2}) = y^k$. If $X_i = 0$, then $y = 1$,

$$\tilde{\pi}_i^0 - \pi_i^0 \leq 0.$$

Otherwise, if $X_i \neq 0$, then $y > 1$, and

$$\tilde{\pi}_i^0 - \pi_i^0 \leq C\frac{y-1}{y^k} = C\frac{y-1}{(y-1+1)^k} \leq C\frac{y-1}{k(y-1)} < C\epsilon. \tag{A.6}$$

Combine (A.5) and (A.6), then

$$|\tilde{\pi}_i^0 - \pi_i^0| \leq max(\delta, C\epsilon) < \delta. \tag{A.7}$$

Now, assume that $\pi_0$ is also estimated by $\hat{\pi}_0$. Let $A = \frac{\sigma}{\sqrt{\sigma^2+\hat{\tau}^2}}\exp(\frac{\hat{M}X_i^2}{2\sigma^2})$, then

$$|\hat{\pi}_i^0 - \tilde{\pi}_i^0| = |\frac{\hat{\pi}_0}{\hat{\pi}_0 + \hat{\pi}_1 A} - \frac{\pi_0}{\pi_1 + \pi_1 A}| = |\frac{(\hat{\pi}_0 - \pi_0)A}{(\hat{\pi}_0 + \hat{\pi}_1 A)(\pi_0 + \pi_1 A)}|$$

The denominator greater than $\hat{\pi}_0\pi_1 A$ implies that $|\hat{\pi}_i^0 - \tilde{\pi}_i^0| < |\frac{\hat{\pi}_0 - \pi_0}{\hat{\pi}_0\pi_1}|$. Since $\hat{\pi}_0$ is consistent for $\pi_0$, for any $\delta > 0$, there $\exists P_0$ such that $\forall p > P_0$, $|\hat{\pi}_0 - \pi_0| < \delta$, then

$$|\hat{\pi}_i^0 - \tilde{\pi}_i^0| \leq D\delta,$$

where $D$ is a constant that only depends on $\pi_0$. Combining this with (A.7), one can get that

$$|\hat{\pi}_i^0 - \pi_i^0| \leq (1+D)\delta, \text{ for all } i = 1, 2, \cdots, p$$

and completes the proof.

**Proof of the theorem**

According to Zhao and Hwang (2008), $FCR_\pi = E\frac{\sum_{i\in\mathcal{R}}P(\theta_i\notin CI_i^{EB}|X)}{R}I(R > 0)$ where $\mathcal{R}$ is the set of index of parameters that are selected and $R$ is the number of selected parameters, i.e. $R = \#\mathcal{R}$. Similarly as formula (A.3) in the proof of theorem 3.1,

$$P(\theta_i \notin CI_i^{EB}|X)$$
$$= \pi_i^0(X)I(\hat{\pi}_i^0(X) < \hat{k}_2) + (1 - \pi_i^0(X))P(\theta_i \notin CI_i^{EB}|X, I_i = 1)$$

In the empirical Bayes interval (10), there exists a positive correction term $-\hat{M}\log\hat{M}\sigma^2$. Dropping this term results in a short interval which enlarges the non-coverage probability, i.e.

$$P(\theta_i \notin CI_i^{EB}|X) \leq P(|\theta_i - \hat{M}X_i)^2 > \hat{M}\sigma^2 z_{q/2}^2).$$

Consequently,

$$P(\theta_i \notin CI_i^{EB}|X) \leq \pi_i^0(X)I(\hat{\pi}_i^0(X) < \hat{k}_2) + (1-\pi_i^0(X))P((\theta_i-\hat{M}X_i)^2 > \hat{M}\sigma^2 z_{q/2}^2|X, I_i = 1).$$

Rearrange the terms in the above formula, one can simply the conditional non-coverage probability $P(\theta_i \notin CI_i^{EB}|X)$ as

$$\pi_i^0(X)(I(\hat{\pi}_i^0(X) < \hat{k}_2) - q) + \pi_i^0(X)(q - P((\theta_i - \hat{M}X_i)^2 > \hat{M}\sigma^2 z_{q/2}^2|X, I_i = 1))$$

$$+ P((\theta_i - \hat{M}X_i)^2 > \hat{M}\sigma^2 z_{q/2}^2|X, I_i = 1).$$

Let

$$\Delta_1 = \frac{\sum_{i \in A} \pi_i^0(X)(I(\hat{\pi}_i^0(X) < \hat{k}_2) - q)}{R},$$

$$\Delta_2 = \frac{\sum_{i \in A} \pi_i^0(X)(q - P((\theta_i - \hat{M}X_i)^2 > \hat{M}\sigma^2 z_{q/2}^2|X, I_i = 1))}{R},$$

and

$$\Delta_3 = \frac{\sum_{i \in A} P((\theta_i - \hat{M}X_i)^2 > \hat{M}\sigma^2 z_{q/2}^2|X, I_i = 1)}{R},$$

then $FCR_\pi$ can be controlled from above by $E(\Delta_1 + \Delta_2 + \Delta_3)$.

Since $\hat{\pi}_0$ and $\hat{\tau}^2$ are obtained by using the method of moments, Delta method implies that $\hat{\pi}_0 - \hat{\pi} = O_p(\frac{1}{\sqrt{p}})$ and $\hat{\tau}^2 - \tau^2 = O_p(\frac{1}{\sqrt{p}})$.

According to Lemma (A.1), for any $\epsilon > 0$, I can always find sufficiently large $P_0$, such that for any $p > P_0$, $(\hat{\tau}^2 - \tau^2)^2 < \delta/3$ and $(\hat{\pi}_i^0(X) - \hat{\pi}_i^0(X))^2 < \delta/3$. Consequently,

$$E\Delta_1 \leq E\frac{\sum_{i \in A} \pi_i^0(X)(I(\pi_i^0(X) < \hat{k}_2 + \sqrt{\delta/3}) - q)}{R} = f(p, \tau^2, \pi_0, \hat{k}_2 + \sqrt{\delta/3}).$$

Since $(\hat{\tau}^2 - \tau^2)^2 + (\hat{\pi}_i^0(X) - \pi_i^0(X))^2 + (\delta/3)^2 \leq \delta$, therefore according to the property of the function $f$,

$$f(p, \tau^2, \pi_0, \hat{k}_2 + \sqrt{\delta/3}) \leq f(p, \hat{\tau}^2, \hat{\pi}_0, \hat{k}_2) + \epsilon \leq \epsilon,$$

Since $\hat{k}_2$ is simulated as the maximum $k_2$ such that $f(p, \hat{\tau}^2, \hat{\pi}_0, k_2) \leq 0$,

$$E\Delta_1 \leq \epsilon. \tag{A.8}$$

For the second term $\Delta_2$,

$$
\begin{aligned}
|\Delta_2| &\leq \frac{\sum_{i \in A} \pi_i^0(X)|q - P((\theta_i - \hat{M}X_i)^2 > \hat{M}\sigma^2 z_{q/2}^2|X, I_i = 1)|}{R} \\
&\leq \frac{\sum_{i \in A} |q - P((\theta_i - \hat{M}X_i)^2 > \hat{M}\sigma^2 z_{q/2}^2|X, I_i = 1)|}{R}.
\end{aligned}
$$

Taking a close look at the term $P((\theta_i - \hat{M}X_i)^2 > \hat{M}\sigma^2 z_{q/2}^2|X, I_i = 1)$, one knows that $(\theta_i|X_i, I_i = 1) \sim N(MX_i, M\sigma^2)$. Therefore one can replace $\theta_i$ by $MX_i + \sqrt{M}\sigma Z$ where $Z$ is a standard normal random variable which is independent of $X_i$. Consequently,

$$P((\theta_i - \hat{M}X_i)^2 > \hat{M}\sigma^2 z_{q/2}^2|X, I_i = 1) = P(|Z - \frac{(\hat{M} - M)X_i}{\sqrt{M}\sigma}| > \sqrt{\frac{\hat{M}}{M}} z_{q/2}|X). \tag{A.9}$$

Assume that $X_{(p)}$ is the observation that has the largest absolute value, then $0 \leq |\frac{(\hat{M}-M)X_i}{\sqrt{M}\sigma}| \leq |\frac{(\hat{M}-M)X_{(p)}}{\sqrt{M}\sigma}|$. Consequently, for any $i = 1, 2, \cdots, p$, (A.9) falls into the range

$$[P(|Z - \frac{|(\hat{M}-M)X_{(n)}|}{\sqrt{M}\sigma}| \geq \sqrt{\frac{\hat{M}}{M}}z_{q/2}), P(|Z| \geq \sqrt{\frac{\hat{M}}{M}}z_{q/2})]. \qquad (A.10)$$

Let $X_i = \sqrt{\sigma^2 + \tau^2}Z_i$, then $Z_i = \pi_0 N(0, \frac{\sigma^2}{\sigma^2+\tau^2}) + \pi_1 N(0,1)$. Furthermore

$$|\frac{(\hat{M}-M)X_{(p)}}{\sqrt{M}\sigma}| = |\frac{\sigma(\hat{\tau}^2 - \tau^2)}{\tau(\hat{\tau}^2 + \sigma^2)}Z_{(p)}|,$$

As a result, the range (A.10) can be rewritten as

$$[P(|Z - |\frac{\sigma(\hat{\tau}^2 - \tau^2)}{\tau(\hat{\tau}^2 + \sigma^2)}Z_{(p)}|| \geq |\sqrt{\frac{\hat{M}}{M}}z_{q/2}), P(|Z| \geq \sqrt{\frac{\hat{M}}{M}}z_{q/2})]. \qquad (A.11)$$

Since the above range applies for all $i$'s, one knows that

$$|\Delta_2| \leq max(|q - P(|Z - |\frac{\sigma(\hat{\tau}^2 - \tau^2)}{\tau(\hat{\tau}^2 + \sigma^2)}Z_{(p)}|| \geq |\sqrt{\frac{\hat{M}}{M}}z_{q/2})|, |q - P(|Z| \geq \sqrt{\frac{\hat{M}}{M}}z_{q/2})|). \qquad (A.12)$$

Since $\hat{\tau}^2 - \tau^2 = O_p(\frac{1}{\sqrt{p}})$, $Z_{(p)} = O(\sqrt{2\log p})$,

$$|\frac{\sigma(\hat{\tau}^2 - \tau^2)}{\tau(\hat{\tau}^2 + \sigma^2)}Z_{(p)}| = o_p(1). \qquad (A.13)$$

The dominated convergence theorem implies that

$$P(|Z - |\frac{\sigma(\hat{\tau}^2 - \tau^2)}{\tau(\hat{\tau}^2 + \sigma^2)}Z_{(p)}|| \geq |\sqrt{\frac{\hat{M}}{M}}z_{q/2})) \to P(|Z| > z_{q/2}) = q,$$

and

$$P(|Z| \geq \sqrt{\frac{\hat{M}}{M}}z_{q/2}) \to q.$$

Applying the dominated convergence theorem again, one can deduce from (A.12) that

$$\limsup_{p\to\infty} E|\Delta_2| \leq 0. \qquad (A.14)$$

Similar arguments apply to $\Delta_3$ and one can show that

$$\Delta_3 \leq P(|Z - \frac{|(\hat{M}-M)X_{(p)}|}{\sqrt{M}\sigma}| \geq \sqrt{\frac{\hat{M}}{M}}z_{q/2}|X)$$

$$= P(|Z - |\frac{\sigma(\hat{\tau}^2 - \tau^2)}{\tau(\hat{\tau}^2 + \sigma^2)}Z_{(p)}|| \geq \sqrt{\frac{\hat{M}}{M}}z_{q/2}).$$

Dominated convergence theorem and (A.13) implies that

$$\limsup_{p \to \infty} E|\Delta_3| \le \lim_{p \to \infty} EP(|Z| \ge z_{q/2}) = q. \tag{A.15}$$

(A.8), (A.14), and (A.15) imply that

$$\limsup_{p \to \infty} FCR_\pi \le q.$$

**Proof of the proposition 4.1.**

Assume that $X_i \sim \pi_0 N(0, \sigma^2) + (1 - \pi_0) N(0, \tau^2 + \sigma^2)$ and $Y_i \sim \pi_0' N(0, \sigma^2) + (1 - \pi_0') N(0, \tau'^2 + \sigma^2)$ where $i = 1, 2, \cdots, p$. Then

$$
\begin{aligned}
& |f(p, \pi_0, \tau^2, k) - f(p, \pi_0', \tau'^2, k')| \\
= {} & E \frac{\sum_{i=1}^{R} \pi_i^0(X)(I(\pi_i^0(X) < k) - q) - \pi_i'^0(Y)(I(\pi_i'^0(Y) < k') - q)}{R} \\
= {} & E \frac{q \sum_{i=1}^{R} (\pi_i^0(X) - \pi_i'^0(Y)) + \sum_{i=1}^{R} (\pi_i^0(X) I(\pi_i^0(X) < k) - \pi_i'^0(Y) I(\pi_i'^0(Y) < k'))}{R}.
\end{aligned}
$$

Since $R$ goes to $\infty$ as $p \to \infty$, therefore by using the law of large number, the inside function of the above expectation converges to $\Delta = qE(\pi_1^0(X) - \pi_1'^0(Y)) + E(\pi_1^0(X) I(\pi_1^0(X) < k) - \pi_1'^0(Y) I(\pi_1'^0(Y) < k'))$ in probability. Since the integral is a bounded function, it is sufficient to show that $\forall \epsilon > 0$, there exists $\delta$, such that $(k' - k)^2 + (\tau'^2 - \tau^2)^2 + (\pi_0' - \pi_0)^2 < \delta$ implies that $|\Delta| < \epsilon$.

In fact, $E\pi_1^0(X) E(P(\theta_0 = 0 | X)) = P(\theta_0 = 0) = \pi_0$. This implies that

$$qE(\pi_1^0(X) - \pi_1'^0(Y)) = q(\pi_0 - \pi_0'). \tag{A.16}$$

Furthermore, direct calculation shows that

$$E(\pi_1^0(X) I(\pi_1^0(X) < k)) = \int_{\pi_1^0(X) < k} \pi_1^0(X) m(X) dX = \pi_0 \int_{\pi_1^0(X) < k} \frac{1}{\sqrt{2\pi\sigma^2}} \exp(-\frac{x^2}{2\sigma^2}) dx.$$

Since $\{\pi_1^0(X) < k\}$ implies that $|X|^2 > \frac{2\sigma^2}{M}(\log \frac{1-k}{k} + \log \frac{\pi_0}{\pi_1 \sqrt{1-M}})$,

$$
\begin{aligned}
& E(\pi_1^0(X) I(\pi_1^0(X) < k) - \pi_1'^0(Y) I(\pi_1'^0(Y) < k')) \\
= {} & P(|N|^2 > \frac{2\sigma^2}{M}(\log \frac{1-k}{k} + \log \frac{\pi_0}{\pi_1 \sqrt{1-M}})) - P(|N|^2 > \frac{2\sigma^2}{M'}(\log \frac{1-k'}{k'} + \log \frac{\pi_0'}{\pi_1' \sqrt{1-M'}})),
\end{aligned}
$$

where $N$ is a standard normal random variable. When $k, k'$ are close to 1, then $\log \frac{1-k}{k} \to -\infty$, therefore, $|E(\pi_1^0(X) I(\pi_1^0(X) < k) - \pi_1'^0(Y) I(\pi_1'^0(Y) < k'))| = 0$ if $k, k' > \epsilon_1$ where $\epsilon_1 < 1$ is close to 1 sufficiently. Similarly, if $k, k'$ are close to 0, then $\log \frac{1-k}{k} \to \infty$. I can choose sufficiently small $\epsilon_0$, such that when $k, k' < \epsilon_0$,

$$P(|N|^2 > \frac{2\sigma^2}{M}(\log \frac{1-k}{k} + \log \frac{\pi_0}{\pi_1 \sqrt{1-M}})) < \frac{\epsilon}{2}, P(|N|^2 > \frac{2\sigma^2}{M'}(\log \frac{1-k'}{k'} + \log \frac{\pi_0'}{\pi_1' \sqrt{1-M'}})) < \frac{\epsilon}{2}.$$

Consequently, $|E(\pi_1^0(X)I(\pi_1^0(X) < k) - \pi_1'^0(Y)I(\pi_1'^0(Y) < k'))| < \epsilon$ when $k, k'$ are either close to 0 or 1.

Furthermore, assume that $0 < \epsilon_0 < k, k' < \epsilon_1 < 1$, then by the continuity of $E(\pi_1^0(X)I(\pi_1^0(X) < k) - \pi_1'^0(Y)I(\pi_1'^0(Y) < k'))$, there exists a small $\delta < \epsilon$, such that $(k' - k)^2 + (\tau'^2 - \tau^2)^2 + (\pi_0' - \pi_0)^2 < \delta$ implies that $|E(\pi_1^0(X)I(\pi_1^0(X) < k) - \pi_1'^0(Y)I(\pi_1'^0(Y) < k'))| < \epsilon$. Combining this with (A.16), one obtains that $|\Delta| < \epsilon$ when $\delta$ is sufficiently small, which completes the proof.

### A.2. Figures and Graphs

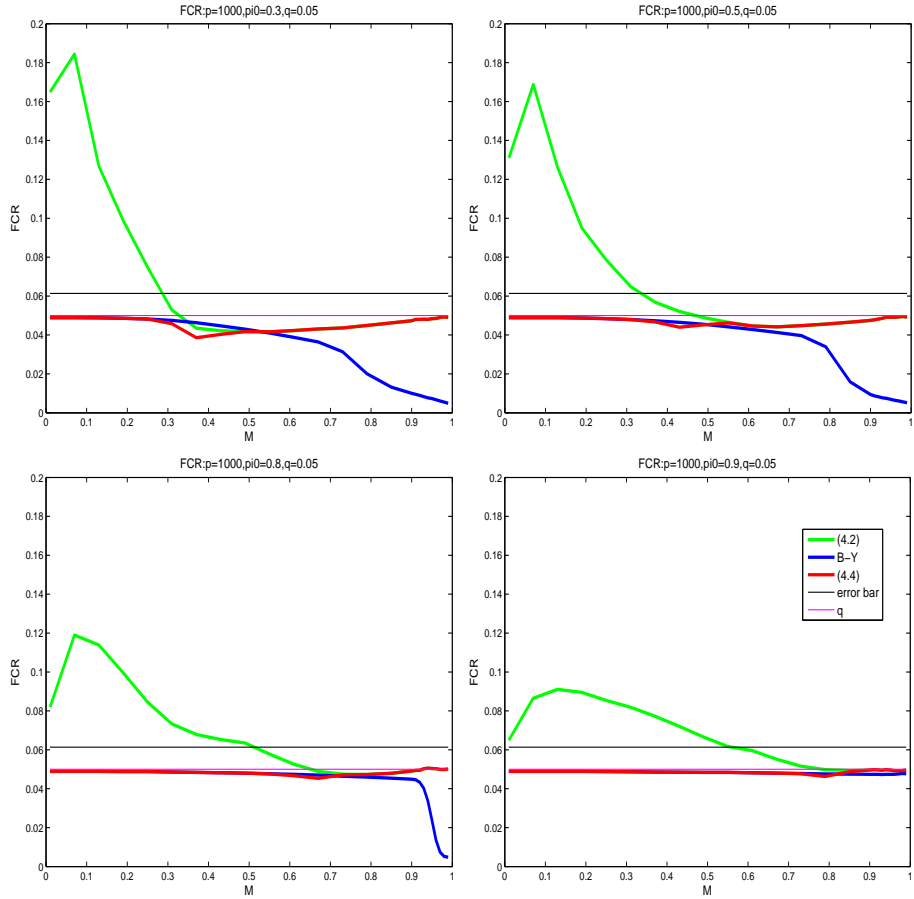

Fig 1. *This figures are the simulated Bayes FCR under different model settings against $M = \frac{\tau^2}{1+\tau^2}$. The dimension is set to be 1000, and top 100 observations after ordering all $X_i$'s according to their magnitude are selected for confidence interval construction. The hyper parameter $\pi_0$ varies among $0.3, 0.5, 0.8$ and $0.9$. The Bayes FCR level that I aim at is 5%. When $\tau^2$ is small, (10) doesn't control the Bayes FCR at 5%; however, the mixed procedure (12) does control the Bayes FCR for any hyper parameters. The portions of the mixture increase as $\pi_0$ increases.*

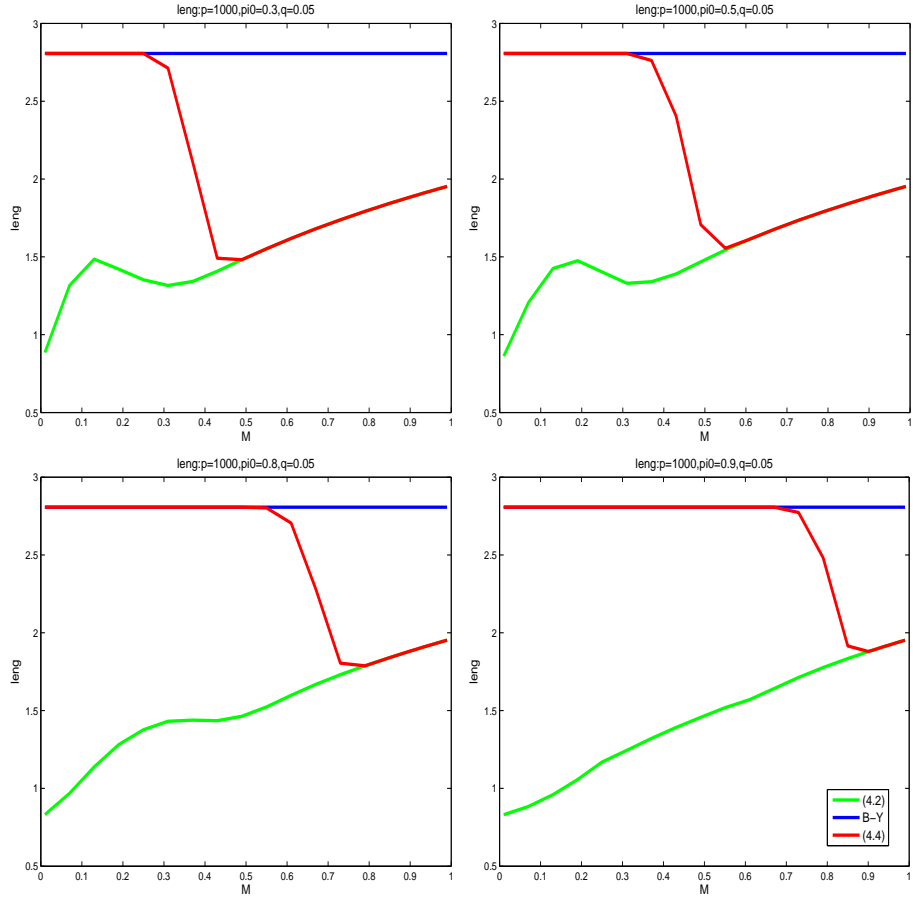

FIG 2. *This figures are the simulated average lengths of different approaches under the same model settings as figure 1. The average length of my procedure is less or equal than B-Y's procedure. In some extremely cases, the average length of (12) is only 54% of that of B-Y's procedure.*

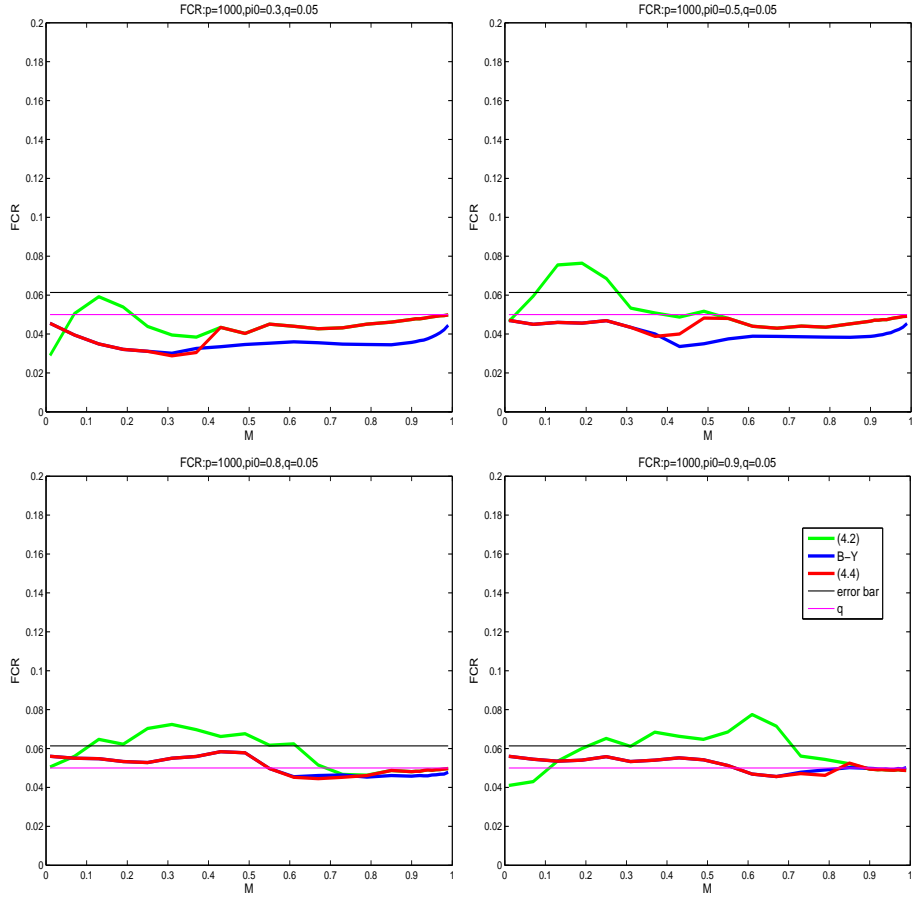

Fig 3. *This figures are the simulated Bayes FCR under different model settings against M =* $\frac{\tau^2}{1+\tau^2}$. *The dimension is set to be 1000. The selection rule is based on Benjamini and Hochberg (1995)'s procedure which aims at controlling the False Discovery Rate to be less or equal than 5%. The hyper parameter $\pi_0$ varies among $0.3, 0.5, 0.8$ and $0.9$. The Bayes FCR level that I aim at is 5% which is represented by the magenta line. When $\tau^2$ is small, (10) doesn't control the Bayes FCR; however, FCRs of the mixed procedure (12) and B-Y's procedure are always less or equal than the error bar which equals to q plus the simulation error.*

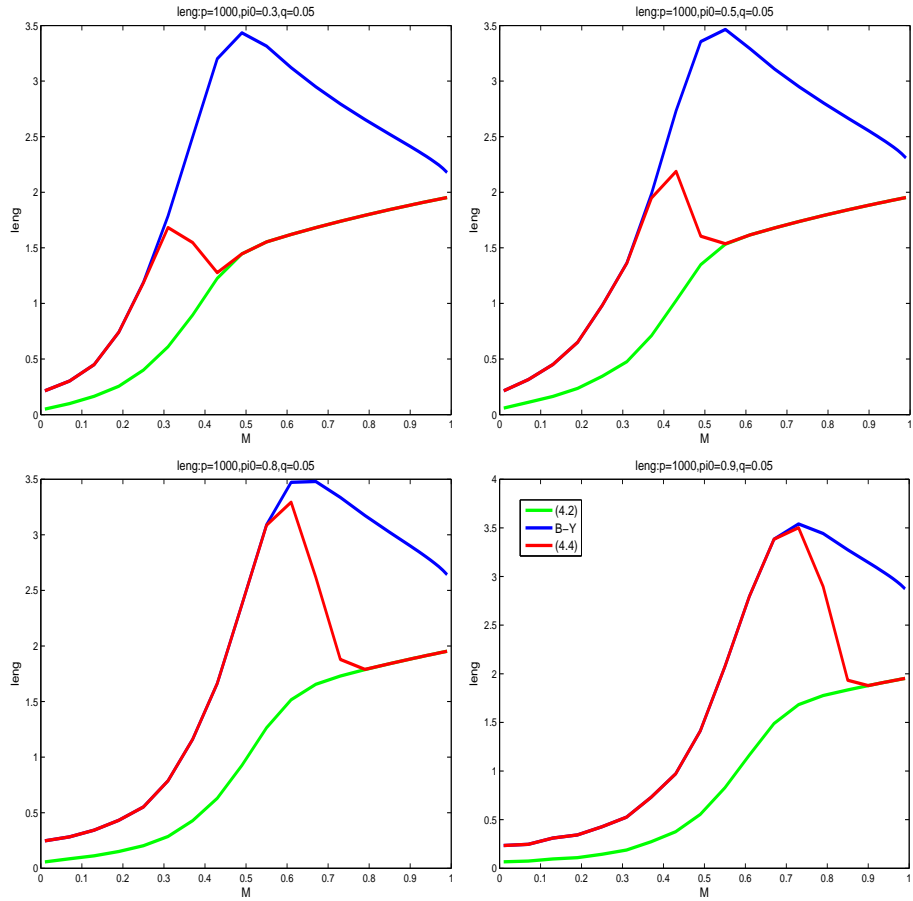

FIG 4. *This figures are the simulated average length of different approaches under the same model as figure 1. The average length of my procedure less than B-Y's procedure. In some extremely cases, the average length of (12) is only 44% of that of B-Y's procedure.*

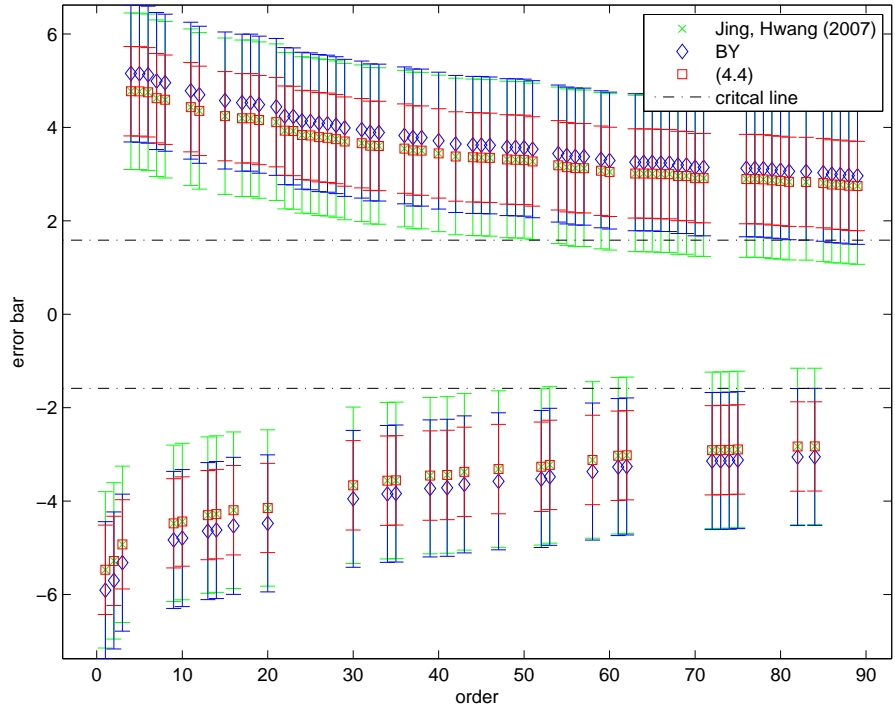

FIG 5. *Three different interval approaches, Qiu and Hwang (2007), B-Y (2005), and (12) are applied to the Synteni data of Kerr, Martin, and Churchill (2000). B-H (1995)'s FDR procedure which aims at finding the genes with differentially expressed levels to be significantly larger or equal than $\log_2 3$ while controlling the False Discovery Rate to be less or equal than 5% is applied to select genes for interval estimation. Among 1285 genes, 89 of them are declared significant and the corresponding intervals are constructed and plotted in this figure. From the figure, one can see that the center of Qiu and Hwang (2007)'s procedure is the same as (12). However, since they aim at controlling the simultaneous coverage coefficient by using Bonferroni's correction, lengths of their intervals are much larger than that of (12). B-Y (2005) centers their intervals at the biased estimator $X_{(i)}$'s. Thus they end up correcting the selection bias by increasing the length and it turned out that their lengths are much larger than that of (12). However, B-Y's length is slightly smaller than Qiu and Hwang (2007)'s procedure.*